\documentclass{article}
\usepackage{latexsym,amssymb,amsmath}

\usepackage{amsfonts}
\newtheorem{theorem}{Theorem}[section]
\newtheorem{corollary}[theorem]{Corollary}
\newtheorem{lemma}[theorem]{Lemma}
\newtheorem{problem}[theorem]{Problem}
\newtheorem{proposition}[theorem]{Proposition}
\newtheorem{rem}[theorem]{Remark}
\newtheorem{example}[theorem]{Example}
\newenvironment{remark}{\begin{rem} \em}{\end{rem}}

\def\Gl{\mathop{\rm Gl}\nolimits}
\def\deg{\mathop{\rm deg }\nolimits}
\def\rank{\mathop{\rm rank}\nolimits}

\newcommand{\se}{\ensuremath{\stackrel{s.e.}{\sim}}}
\newcommand{\s}{\ensuremath{\stackrel{s}{\sim}}}

\newenvironment{matriz}[1]{\left[ \begin{array}{#1}}{\end{array} \right]}

\newcommand{\FF}{\mathbb F}

\newcommand{\ZZ}{\mathbb Z}

\newcommand{\MSC}[1]{\textbf{\textit{MSC}} #1}
\title{Fixed rank perturbations of regular matrix pencils}
\author{Itziar Baraga\~na\footnote{Departamento de Ciencia de la
Computaci\'on e IA,
Universidad del Pa\'{\i}s Vasco, UPV/EHU,
Apdo. 649,
20080 Donostia-San Sebasti\'an, Spain, e-mail: itziar.baragana@ehu.es.
Partially supported by MINECO:
MTM2017-83624-P,
MTM2017-90682-REDT, and UPV/EHU: GIU16/42.
},
Alicia Roca\footnote{Corresponding author.
Departamento de  Matem\'atica Aplicada,  IMM,
 Universitat Polit\`ecnica de Val\`encia,  46022 Valencia, Spain,
e-mail: aroca@mat.upv.es.
Partially supported by MINECO:
MTM2017-83624-P,
MTM2017-90682-REDT.
}
}

\date{}

\begin{document}

\maketitle

\begin{abstract}
A characterization of  the structure of a regular matrix pencil obtained by a bounded  rank  perturbation of another regular matrix pencil has been recently obtained. The result generalizes the solution for the  bounded rank perturbation problem of a square constant  matrix. 
When comparing the fixed rank perturbation problem of a constant matrix with the bounded rank perturbation problem we realize that both problems are of different nature;  the first one is more restrictive.  In this paper we  characterize the  structure of a regular matrix pencil obtained by a fixed rank  perturbation of another regular matrix pencil.
We apply the result to find necessary and sufficient conditions for the existence of a fixed rank perturbation such that the perturbed pencil has a prescribed determinant. The results hold over   fields with sufficient number of elements.\end{abstract}


\MSC 15A22, 47A55, 15A18.

\section{Introduction}
\label{secintroduction}

Low rank perturbations of matrix pencils have been widely studied, and  the problem has recently deserved   the attention of several authors, as we will see in the next references. Given a matrix pencil  $A(s)$ and a nonnegative integer $r$, the problem consists in  characterizing the Kronecker structure of  $A(s)+P(s)$, where $P(s)$ is a matrix pencil of bounded ($\rank (P(s))\leq r$) or fixed rank ($\rank (P(s))= r$). 

Some authors focus their research on {\em generic} perturbations; it means that the perturbation pencil $P(s)$ belongs  to an open and dense subset of the set of pencils of bounded or fixed rank  (for this approach see for instance \cite{TeDo07, TeDo16,  TeDoMo08, MoDo03, Sa02, Sa04} and the references therein).

In other papers   the pencil $P(s)$ is allowed to be an  {\em arbitrary} perturbation belonging to the whole set of pencils of  bounded or fixed rank.
Within this framework and for bounded rank perturbations,  the problem has  been solved in \cite{Silva88_1, Za91} for pencils $A(s)=sI-A$, $P(s)=P$, with $A, P$ constant matrices (see Proposition \ref{propSiZa} and Corollary \ref{corless} below). In the  earlier work  \cite{Th80},  the same problem was solved for $r=1$.  A solution for quasi-regular matrix pencils of the form $A(s)=\begin{matriz}{cc} sI_n-A_1 &  A_2 \end{matriz}$ and constant perturbation pencils $P(s)=\begin{matriz}{cc} P_1 & P_2 \end{matriz}$  has been obtained in \cite{DoSt14}. For regular pencils $A(s)$  and $A(s)+P(s)$, a solution to  the problem has recently been given in \cite{BaRo18} (see Proposition \ref{propless}) (see also \cite{BaRo18} for further references on the problem). 

Concerning fixed rank perturbations, the problem has been solved in \cite{Silva88_1}  when $A(s)=sI-A$ and $P(s)=P$ is a constant matrix. The  result obtained holds over algebraically closed fields (see Proposition \ref{propSi2} below). When comparing the characterization of the solutions of the bounded (\cite{Silva88_1, Za91}) and fixed rank (\cite{Silva88_1}) perturbation problems, we observe that an extra condition appears in the fixed rank case, which proves that the two problems are of different nature.

In this paper we deal with regular matrix pencils and we require that  $P(s)$ is a matrix pencil of fixed rank.  More precisely, the first problem we solve is the following:

\begin{problem}\label{fixproblem}
Given two regular matrix pencils  $A(s), B(s)\in \FF[s]^{n \times n}$
and a nonnegative integer $r$, $r\leq n$,
 find necessary and sufficient conditions for the existence of a
matrix pencil
 $P(s)\in \FF[s]^{n \times n}$ such that $\rank (P(s))= r$ and 
 $A(s)+P(s)$ is strictly equivalent to $B(s)$.
\end{problem}

Recall that when $A(s)$ is a regular matrix pencil, the Kronecker structure of $A(s)$  is formed by its homogeneous invariant factors, and it is known as the Weierstrass structure of the pencil (see Theorem \ref{theoWei}). 

A solution to Problem \ref{fixproblem} is given in Theorems \ref{theon1} and \ref{maintheorem}. 
Unlike what happens when perturbing pencils of the form $sI-A$ with constant matrices, 
in this case the solutions to the bounded and  fixed rank perturbation problems are characterized in terms of the same conditions. This is due to the fact that as the perturbation matrix can be a matrix pencil, it introduces some more freedom that in the constant perturbation problem. But, the fact of being a more restrictive problem determines extra needs for achieving a solution, and in this case proofs are  more demanding. To solve it under the same conditions of the bounded case, some specific technical lemmas must be introduced; nothing similar was required in the bounded case.

\medskip

The solution  to Problem \ref{fixproblem} obtained allows us to solve the following eigenvalue placement problem:

\begin{problem}\label{prfixdet}
  Given a regular matrix pencil  $A(s)\in \FF[s]^{n \times n}$, a nonnegative integer $r$, $r\leq n$, and a monic polynomial  $0\neq q(s)\in \FF[s]$ with $\deg(q(s))\leq n$,
 find necessary and sufficient conditions for the existence of a
 matrix pencil
 $P(s)\in \FF[s]^{n \times n}$ such that $\rank (P(s))= r$ and
$\det(A(s)+P(s))=kq(s)$, with $k\in \FF$.
\end{problem}

A solution to Problem \ref{prfixdet} is given in Theorem \ref{thomainplacement} (see also Corollary   \ref{corplacement}
and Remark \ref{rempln1}). 
An analogous problem was solved in   \cite{BaRo18}  in the case that $\rank (P(s))\leq r$. For $r=1$, see also \cite{GeTr17}.

\medskip
The paper is organized as follows. In Section \ref{secpreliminaries} we introduce the notation, basic definitions and preliminary results. In Section \ref{secmain} we solve  Problem \ref{fixproblem}, first for pencils not having infinite elementary divisors and then for the general case. A solution to Problem \ref{prfixdet} is given in Section \ref{secplacement}. Finally, in Section \ref{conclusions} we summarize the main contributions of the paper.

\section{Notation and preliminary results}
\label{secpreliminaries}

The section contains three subsections, where we introduce notation and basic definitions (Subsection \ref{subsecnotation}), some results concerning matrix pencils (Subsection \ref{subsecpencils}), and previous results about matrix or pencil perturbations of bounded or fixed rank (Subsection \ref{subseclowrank}).

\subsection{Notation and basic definitions}
\label{subsecnotation}

Let $\FF$ be a field. $\FF[s]$ denotes the ring of polynomials in the indeterminate $s$ with coefficients in $\FF$,  $\FF[s, t]$  the ring of polynomials in two indeterminates $s, t$ with coefficients in $\FF$,
and  $\FF^{m\times n}$,  $\FF[s]^{m\times n}$ and $\FF[s, t]^{m\times n}$ the vector spaces  of $m\times n$ matrices with elements in $\FF$, $\FF[s]$ and $\FF[s, t]$, respectively.
$\Gl_n(\FF)$ is the general linear group of invertible matrices in $\FF^{n\times n}$.

\medskip

The number of elements of a finite set $I$ will be denoted by $\mid I \mid$.
If $G$ is a matrix in $\FF^{m\times n}$, $I\subseteq \{1, \dots, m\}$, and
$J\subseteq \{1, \dots, n\}$, with $\mid I \mid=r$ and  $\mid J \mid=s$, then
$G(I, J)$ denotes the $r\times s$ submatrix of  $G$ formed by the rows in $I$ and the columns in $J$.
Similarly, $G(I, :)$ is the $r\times n$ submatrix of $G$ formed by the rows in $I$ and $G(:, J)$
is  the $m \times s$ submatrix of $G$ formed by the columns in $J$.

If $\mid I \mid=\mid J \mid$ and $\det (G(I, J))\neq 0$, then the {\em Schur complement of $G(I, J)$ in $G$} is
$$ G/G(I,J)=
    G(I^c, J^c)-G(I^c, J)G(I, J)^{-1}G(I, J^c),
    $$
    where $I^c=\{1, \dots, m\}\setminus I$ and $J^c=\{1, \dots, n\}\setminus J$ (see \cite{Ando87}).
    It is satisfied that
    $$\rank (G)=\rank (G(I, J)) + \rank (G/G(I, J)),$$ and if $m=n$, 
    $$\det (G)=\pm \det (G(I, J)) \det(G/G(I, J)).$$

    \medskip
Given a polynomial matrix $G(s)\in \FF[s]^{m\times n}$, the {\em degree} of $G(s)$, denoted by $\deg (G(s))$, is the  maximum  of  the degrees of its entries. The {\em normal rank} of $G(s)$, denoted by $\rank (G(s))$,  is the order of the largest non identically zero minor of $G(s)$, i.e. it is the rank of $G(s)$ considered as a matrix on the field of fractions of $\FF[s]$. If $\rank (G(s))=\rho$, the {\em determinantal divisor of order $k$} of  $G(s)$, denoted by $D_k(s)$, is the monic greatest common divisor  of the minors of order  $k$ of $G(s)$,  $1\leq k \leq \rho$. The determinantal divisors satisfy 
$D_{k-1}(s)\mid D_k(s)$, $1\leq k \leq \rho$ ($D_0(s):=1$) and  the {\em invariant factors of $G(s)$} are the monic polynomials
$$
\gamma_k(s)=\frac{D_k(s)}{D_{k-1}(s)}, \quad 1\leq k \leq \rho.
$$
We will take $\gamma_i(s):=1$ for $i<1$ and  $\gamma_i(s):=0$ for $i>\rho$.

\medskip
A matrix    $U(s)\in \FF[s]^{n\times n}$ is  {\em unimodular} if $0\neq \det (U(s))\in \FF$, equivalently $U(s)$ is a unit in the ring  $\FF[s]^{n\times n}$.
Two polynomial matrices $G(s), H(s)\in \FF[s]^{m\times n}$ are {\em equivalent} ($G(s)\sim H(s)$) if there exist unimodular matrices $U(s)\in \FF[s]^{m\times m}$, $ V(s) \in \FF[s]^{n\times n}$ such that
$G(s)=U(s)H(s)V(s)$. A complete system of invariants for the equivalence of polynomial matrices is formed by the invariant factors, i.e. two polynomial matrices $G(s), H(s)\in \FF[s]^{m\times n}$ are equivalent if and only if they have the same invariant factors.

\medskip
Given a square matrix $G\in \FF^{n \times n}$, the {\em invariant factors of $G$} are the invariant factors of the polynomial matrix $sI_n-G$.
Two square matrices $G, H\in \FF^{n\times n}$ are {\em similar} ($G\s H$) if there exists an invertible matrix $Q\in \Gl_n(\FF)$,  such that
$G=QHQ^{-1}$. It is well known that $G\s H$ if and only if $sI_n-G\sim sI_n-H$, i.e. if and only if $G$ and $H$ have the same invariant factors (see, for instance, \cite[Ch. 6, Theorem 7]{Ga74}).

\subsection{Matrix pencils}
\label{subsecpencils}

We review now some basic definitions and results about matrix pencils.
For details see,  for example, \cite[Ch. 12]{Ga74}.

\medskip
A {\em matrix pencil} is a  polynomial matrix $G(s)\in \FF[s]^{m\times n}$ with $\deg(G(s))\leq1$. The pencil is  {\em regular} if $m=n$ and $\det (G(s))$ is a  non zero polynomial. Otherwise it is  {\em singular}.

Two matrix pencils
$G(s)=G_0+sG_1, H(s)=H_0+sH_1\in \FF[s]^{m\times n}$ are {\em strictly equivalent} ($G(s)\se H(s)$) if there exist invertible matrices $Q\in \Gl_m(\FF)$,   $R\in \Gl_n(\FF)$ such that
$G(s)=QH(s)R$.

It is immediate that if $G(s)\se H(s)$ then $G(s)\sim H(s)$.
Moreover, if $n=m$, $\det (G_1)\neq 0$ and $\det (H_1 )\neq 0$, then $G(s)\se H(s)$ if and only if
$G(s)\sim H(s)$
(see, for instance, \cite[Ch.12, Theorem 1]{Ga74}).

\medskip

Given  $G(s)=G_0+sG_1\in \FF[s]^{m\times n}$,  with $\rho = \rank(G(s))$, the {\em homogeneous pencil associated to $G(s)$} is
$$G(s, t)=tG_0+sG_1\in \FF[s, t]^{m\times n},$$
and the {\em homogeneous determinantal divisor of order $k$} of  $G(s)$,  denoted by $\Delta_k(s, t)$, is the  greatest common divisor  of  the minors of order  $k$ of $G(s, t)$,
$1\leq k \leq \rho$. We will assume that $\Delta_k(s, t)$ is monic with respect to $s$.
The homogeneous determinantal divisors of $G(s)$ are homogeneous polynomials and  $\Delta_{k-1}(s, t)\mid \Delta_k(s,t)$, $1\leq k \leq \rho$ ($\Delta_0(s,t):=1$).
The {\em  homogeneous invariant factors of $G(s)$}  are the homogeneous polynomials
$$
\Gamma_k(s, t)=\frac{\Delta_k(s,t)}{\Delta_{k-1}(s, t)}, \quad 1\leq k \leq \rho.
$$

If $\gamma_1(s)\mid \dots\mid \gamma_\rho(s)$  are the invariant factors  $G(s)$,
then
$$
\gamma_i(s)=\Gamma_i(s,1), \quad 1\leq i \leq \rho, $$
and
$$ 
\Gamma_i(s,t)=t^{m_i(\infty, G(s))}t^{\deg(\gamma_i)}\gamma_i(\frac{s}{t}),\quad 1\leq i \leq \rho,
$$
for some integers $0\leq m_1(\infty, G(s)) \leq \dots\leq m_\rho(\infty, G(s)) $.
Hence  $\Gamma_1(s,t)\mid \dots \mid\Gamma_{\rho}(s, t)$.
We take $\Gamma_i(s,t):=1$ for $i<1$ and
$\Gamma_i(s,t):=0$ for $i>\rho$.

If  $m_i(\infty, G(s))>0$, then  $t^{m_i(\infty, G(s))}$ is an {\em infinite elementary divisor} of $G(s)$.
The infinite elementary divisors of $G(s)$  exist if and only if $\rank (G_1)< \rank (G(s))$.

We denote by $\overline{ \FF}$ the algebraic closure of $\FF$.
The {\em spectrum} of
$G(s)=G_0+sG_1\in \FF[s]^{m\times n}$ is defined as
$$
\Lambda(G(s))=\{\lambda\in \overline{\FF}\cup\{\infty\}\; : \; \rank (G(\lambda))< \rank (G(s))\},
$$
where we agree that $G(\infty)=G_1$. The elements $\lambda\in \Lambda(G(s))$ are the {\em eigenvalues} of $G(s)$.

The invariant factors and the homogeneous invariant factors of  $G(s)$ can be written as
\begin{equation} \label{gammas}
\gamma_i(s)=\prod_{\lambda\in \Lambda(G(s))\setminus\{\infty\}}(s-\lambda)^{m_i(\lambda, G(s))}, \quad 1\leq i \leq \rho,
\end{equation}
and
\begin{equation} \label{Gammas}
\Gamma_i(s, t)=t^{m_i(\infty, G(s))}\prod_{\lambda\in \Lambda(G(s))\setminus\{\infty\}}(s-\lambda t)^{m_i(\lambda, G(s))}, \quad 1\leq i \leq \rho.
\end{equation}
For $\lambda \in \Lambda(G(s))$, the integers
$0\leq m_1(\lambda, G(s))\leq \dots\leq   m_\rho(\lambda, G(s))$ are called the
{\em partial multiplicities at $\lambda$ of $G(s)$}.
If  $\lambda \in \overline{\FF}\setminus \Lambda(G(s))$,
we put $m_1(\lambda, G(s))=\dots=m_\rho(\lambda,  G(s))=0$.
For $\lambda \in \overline{\FF}\cup \{\infty\}$, we  will agree that $m_i(\lambda, G(s))=0$ for $i<1$ and  $m_i(\lambda, G(s))=\infty$ for $i>\rho$.

\medskip

When a matrix pencil has infinite elementary divisors, we can perform a change of variable which turn it into a new pencil without infinite structure. This will be done in Section \ref{secmain}, and we will need the following results, which can be found in \cite{CaSil91}.

Let
$
X=
\begin{bmatrix}
x&y\\
z&w
\end{bmatrix} \in \Gl_2(\FF)
$.
For a matrix pencil $G(s)=sG_1+G_0\in \FF[s]^{m \times n}$ and  an homogeneous polynomial $\Phi(s,t)\in \FF[s, t]$
we  define:
$$P_X(sG_1+G_0)=s(xG_1+zG_0)+(yG_1+wG_0)\in \FF[s]^{m \times n},$$
$$\Pi_X(\Phi)(s,t)=\Phi (sx+ty, sz+tw)\in \FF[s, t].$$

\begin{lemma}
[\mbox{\cite[Lemma 6]{CaSil91}}]\label{lemmaCasi6}
The functions $P_X, 
\Pi_X$ are invertible and
$$\left(P_X\right)^{-1}=P_{X^{-1}}, \quad 
\left(\Pi_X\right)^{-1}=\Pi_{X^{-1}}.$$
\end{lemma}

\begin{lemma}
[\mbox{\cite[Lemma 7]{CaSil91}}]\label{lemmaCasi7}
Let $\Phi(s,t), \Psi(s,t)\in \FF[s, t]$ be homogeneous polynomials.
Then,
$\Phi(s,t)\mid \Psi(s,t)$ if and only if $\Pi_X(\Phi)(s,t)\mid\Pi_X( \Psi)(s,t)$.
\end{lemma}

\begin{lemma}[\mbox{\cite[Lemma 9]{CaSil91}}]\label{lemmaCasi9}
Let
$G(s)=sG_1+G_0, H(s)=sH_1+H_0\in \FF[s]^{m \times n}$. Then $G(s)\se H(s)$ if and only if $P_X(G(s))\se P_X(H(s)).$

\end{lemma}

\begin{lemma}[\mbox{\cite[Lemma 10]{CaSil91}}]\label{lemmaCasi10}
Let
$G(s)=sG_1+G_0\in \FF[s]^{m \times n}$, $\rho=\rank (G(s))$. Let $\Gamma_1(s,t) \mid \ldots \mid \Gamma_{\rho}(s,t)$ be the homogeneous invariant factors of $G(s)$. Then the  homogeneous invariant factors of
$P_X(G(s))$ are $\Pi_X(\Gamma_1)(s,t) \mid \ldots \mid \Pi_X( \Gamma_{\rho})(s,t)$.

\end{lemma}

\begin{remark}\label{remrank}
Observe that
\begin{enumerate}
\item[(i)] $\rank (P_X(G(s)))=\rank (G(s)).$
\item[(ii)]
In   Lemma \ref{lemmaCasi10},
 $\Pi_X(\Gamma_i)(s,t)$ are not necessarily monic with respect to $s$. In fact,  $\Pi_X(\Gamma_i)(s,t)$ are the homogeneous invariant factors of
$P_X(G(s))$ multiplied by a constant $0\neq k_i\in \FF$.
\end{enumerate}
\end{remark}

\medskip

In this paper we deal with regular matrix pencils. The following theorem states that
the homogeneous invariant  factors  form a complete system of invariants for the strict equivalence of regular pencils. A proof can be found in
\cite[Ch. 12]{Ga74} for infinite fields and in \cite[Ch. 2]{Ro03} for arbitrary fields.
\begin{theorem}[\mbox{Weierstrass}]\label{theoWei}
  Two regular matrix pencils are strictly equivalent if and only if they have the same homogeneous invariant factors.
\end{theorem}

\medskip

For regular matrix pencils, expressions (\ref{gammas}) and (\ref{Gammas}) allow us to write 
$$
\det (G(s))=\prod_{i=1}^n\gamma_i(s)=
\prod_{\lambda\in \Lambda(G(s))\setminus\{\infty\}}(s-\lambda)^{\mu_a(\lambda, G(s))},
$$
$$
\det (G(s, t))
=\prod_{i=1}^n\Gamma_i(s, t)=
t^{\mu_a(\infty, G(s))}\prod_{\lambda\in \Lambda(G(s))\setminus\{\infty\}}(s-\lambda t)^{\mu_a(\lambda, G(s))},
$$
where, for $\lambda \in \overline{\FF}\cup \{\infty\}$,
$\mu_a(\lambda, G(s))= \sum _{i=1}^n m_i(\lambda, G(s))$ is
the {\em algebraic multiplicity} of $\lambda$ in  $G(s)$.
Notice that $\deg(\det (G(s, t)))=n$ and $\deg(\det (G(s)))=n-\mu_a(\infty, G(s)).$

\medskip

Finally, given an homogeneous polynomial $\Gamma(s, t)$,
we will use the following notation
 $$\Lambda(\Gamma(s,t)):=\{\lambda \in \overline{\FF}\cup \{\infty\}\; : \; 
\Gamma(\lambda,1)=0\},$$
where $\Gamma(\infty, 1):=\Gamma(1, 0)$.
With this notation, if $G(s)\in \FF[s]^{n\times n}$ is a regular matrix pencil with  $\Gamma_1(s,t)\mid \dots \mid\Gamma_{n}(s, t)$ homogeneous invariant factors, then
$$
\Lambda(G(s))=\Lambda(\Gamma_n(s,t))=\Lambda(\Gamma_1(s,t) \dots \Gamma_{n}(s, t)).
$$

Also, for a polynomial $q(s)\in \FF[s]$ with $\deg(q(s))\leq n$, we define
$$\Lambda^n(q(s)):=\{\lambda \in \overline{\FF}\; : \; q(\lambda)=0\}\mbox{  if } \deg(q(s))=n,$$
$$\Lambda^n(q(s)):=\{\lambda \in \overline{\FF}\; : \; q(\lambda)=0\}\cup \{\infty\}\mbox{  if } \deg(q(s))<n.$$

\subsection{Rank perturbations of square matrices and regular matrix pencils}
\label{subseclowrank}

The problem of characterizing the Weierstrass structure of a regular matrix pencil obtained by a bounded rank perturbation of another regular matrix pencil
(i.e. Problem \ref{fixproblem} with the relaxed condition $\rank (P(s))\leq r$) was solved in \cite{BaRo18}. The key point in the obtention of the solution was the next result. It was proven in   \cite{Za91} and in \cite{Silva88_1} under another formulation. We present here the version of  \cite{Za91}.

\begin{proposition}
[\mbox{\cite[Theorem 1]{Silva88_1}, \cite[Theorem 3]{Za91}}]
\label{propSiZa}
Let $A, B \in \FF^{n\times n}$ and let $\alpha_1(s)\mid \dots \mid \alpha_n(s)$ and $\beta_1(s)\mid \dots \mid \beta_n(s)$ be the  invariant factors of $A$ and $B$, respectively.  Let $r$ be a nonnegative integer.
Then there exists a matrix
$P\in \FF^{n\times n}$ such that  $\rank (P)\leq r$  and $A+P$ has $\beta_1(s)\mid \dots \mid \beta_n(s)$ as invariant factors if and only if
\begin{equation}\label{eqintif}
\beta_{i-r}(s)\mid \alpha_i(s)\mid\beta_{i+r}(s), \quad 1\leq i \leq n.
\end{equation}
\end{proposition}

Bearing in mind that
$$A+ P\s B\Leftrightarrow sI_n+A+ P\sim sI_n+B\Leftrightarrow sI_n+A+ P\se sI_n+B,$$
we obtain the following corollary.

\begin{corollary}
\label{corless}
 Let $A(s)=sI_n+A, B(s)=sI_n+B\in \FF[s]^{n \times n}$.
Let $\alpha_1(s)\mid \dots \mid \alpha_n(s)$ and
$\beta_1(s)\mid \dots \mid \beta_n(s)$
 be the invariant factors of
$A(s)$ and $B(s)$, respectively.
 Let $r$ be a nonnegative integer.
Then there exists a matrix
$P\in \FF^{n\times n}$ such that  $\rank (P)\leq r$  and $A(s)+P\se B(s) $
if and only if (\ref{eqintif}) holds.
\end{corollary}
The next proposition is the generalization of  Proposition \ref{propSiZa} to regular matrix pencils obtained in \cite{BaRo18}.

\begin{proposition}
  [\mbox{\cite[Theorem 4.12]{BaRo18}}]
  \label{propless}
 Let $A(s), B(s)\in \FF[s]^{n \times n}$ be regular matrix pencils.
Let $\phi_1(s,t)\mid \dots\mid \phi_n(s, t)$ and $\psi_1(s,t)\mid \dots\mid \psi_n(s, t)$ be the homogeneous invariant factors of
$A(s)$ and $B(s)$, respectively, and assume that  $ \FF\cup \{\infty\}\not \subseteq \Lambda(A(s))\cup \Lambda(B(s))  $. Let $r$ be a nonnegative integer.
There exists a matrix pencil  $P(s)\in \FF[s]^{n \times n}$ such that $\rank (P(s))\leq r$ and $A(s)+P(s)\se B(s)$ if and only if
\begin{equation}\label{interlacinghomogr1}
\phi_{i-r}(s, t)\mid \psi_i(s, t)\mid\phi_{i+r}(s, t), \quad 1\leq i \leq n.
\end{equation}
\end{proposition}

From this proposition  we can derive the following  result.

\begin{corollary}\label{corminrbounded2}
 Let $A(s), B(s)\in \FF[s]^{n \times n}$ be regular matrix pencils.
Let $\phi_1(s,t)\mid \dots\mid \phi_n(s, t)$ and $\psi_1(s,t)\mid \dots\mid \psi_n(s, t)$ be the homogeneous invariant factors of
$A(s)$ and $B(s)$, respectively, and assume that  $ \FF\cup \{\infty\}\not \subseteq \Lambda(A(s))\cup \Lambda(B(s))  $. 
Let
$$
r_0=\min\{r\geq 0 \; : \;  \phi_{i-r}(s,t)\mid\psi_{i}(s,t)\mid\phi_{i+r}(s,t), \quad 1\leq i \leq n\}.
$$
Then there exists a matrix pencil  $P(s)\in \FF[s]^{n \times n}$ such that $\rank (P(s))= r_0$ and $A(s)+P(s)\se B(s)$.
\end{corollary}

In this paper we will show that for any $r$ , $r_0\leq r\leq n$,  there exists a matrix pencil  $P(s)\in \FF[s]^{n \times n}$ such that $\rank (P(s))= r$ and $A(s)+P(s)\se B(s)$ (see Corollary \ref{corrminfixed2}).

\medskip
When $\FF$ is an algebraically closed field
the possible similarity class  of a square matrix obtained by a fixed rank perturbation of another square matrix was characterized in \cite{Silva88_1}. The result is presented in the next proposition;  the statement is different from the original one and more adapted to our problem.

\begin{proposition} 
[\mbox{\cite[Theorem 2]{Silva88_1}}]
\label{propSi2}
Suppose that $\FF$ is  algebraically closed.
Let $A, B \in \FF^{n\times n}$ and let $\alpha_1(s)\mid \dots \mid \alpha_n(s)$ and $\beta_1(s)\mid \dots \mid \beta_n(s)$ be the  invariant factors of $A$ and $B$, respectively.  Let $r$ be a nonnegative integer, $r\leq n$.
Then there exists a matrix
$P\in \FF^{n\times n}$ with   $\rank (P)= r$  such that $A+P$ has $\beta_1(s)\mid \dots \mid \beta_n(s)$ as invariant factors if and only if
(\ref{eqintif}) is satisfied and
\begin{equation}\label{eq2Si}
r\leq \min \{\rank (A-\lambda I_n)+\rank (B-\lambda I_n)\; : \; \lambda \in \FF\}.
\end{equation}
\end{proposition}

\medskip
As mentioned  in the Introduction section, the aim of this paper is to solve an analogous problem to that solved in Proposition \ref{propSi2}  for regular matrix pencils.
When $\FF$ is  algebraically closed, if $A(s)=sI_n+A$ and $B(s)=sI_n+B$, by Proposition \ref{propSi2} conditions  (\ref{interlacinghomogr1}) and
(\ref{eq2Si}) are sufficient for the existence of a matrix pencil  $P(s)\in \FF[s]^{n\times n}$ such that  $\rank (P(s))= r$  and $A(s)+P(s)\se B(s)$.
Nevertheless, (\ref{eq2Si}) is not a necessary condition, as we can see in the next example.
\begin{example}
  Let $c\in \FF$ ($\FF$ algebraically closed),  $A=B=cI_n$,  $r$ an integer, $0<r\leq n$ and $P(s)=\begin{bmatrix}I_r&0\\0&0\end{bmatrix}(sI_n+A)$. Then
  $\rank (P(s))=r$ and
  $$
sI_n+A+P(s)=
\begin{bmatrix}2I_r&0\\0&I_{n-r}\end{bmatrix}
(sI_n+A)\se sI_n+A=sI_n+B,$$
but
   $$\min\{\rank(A-\lambda I_n)+ \rank(B-\lambda I_n)\; : \; \lambda \in \FF\}=0<r.$$

  \end{example}

\section{Fixed rank perturbation for regular matrix pencils}
\label{secmain}

In this section we  give a complete solution to  Problem \ref{fixproblem}  under the same restriction on the field $\FF$ as in Proposition \ref{propless}.
According to  this proposition, the interlacing conditions (\ref{interlacinghomogr1}) are necessary. We prove  that they are also sufficient, except when $\FF$ is a finite field with $\mid \FF \mid=2$ and $r=n=1$.

Following the strategy of \cite{BaRo18}, we start analyzing the case when the pencils $A(s)$, $B(s)$ do not have infinite elementary divisors.

\subsection{Pencils $A(s)$, $B(s)$ without infinite elementary divisors}
\label{subsecsiasib}

First, we analyze the case when $r=n$, then when $r<n$.
Observe that conditions (\ref{interlacinghomogr1}) are trivially fulfilled for $r=n$.
We  prove in Proposition \ref{proprn}  that for  regular pencils $A(s)=sI_n+A, B(s)=sI_n+B\in \FF[s]^{n\times n}$, $n\geq 2$, there always exists a regular pencil $P(s)\in \FF[s]^{n\times n}$ such that $A(s)+P(s)\se B(s)$. In order to do that  we need the following technical lemma.

\begin{lemma}\label{lemmaei+e}
Let $n\geq 2$. Then there exists a matrix
$E_n\in \Gl_n(\FF)$ such that $I_n+E_n\in \Gl_n(\FF)$.
\end{lemma}

{\bf Proof.}
We prove the result by induction on $n$.

If $n=2$, put $E_2=\begin{bmatrix}1&1\\1&0\end{bmatrix}$. Then
 $E_2, I_2+E_2\in \Gl_2(\FF)$.

Assume that  there exists  $E_p\in \Gl_p(\FF)$ such that $I_p+E_p\in \Gl_p(\FF)$ and let $n=p+1$.

Obviously,  $E_p\neq I_p+E_p$. Therefore, if $R=E_p^{-1}-(I_p+E_p)^{-1}$, then
$R\neq 0$. Let $i, j\in \{1, \dots, p\}$ be such that $R(i,j)\neq 0$ and let $w=-1+e_i^tE_p^{-1}e_j\in \FF$.
We define
$$
E_{p+1}=\begin{bmatrix}
E_p&e_j\\
e_i^t&w
\end{bmatrix}\in \FF^{(p+1)\times (p+1)}.
$$
Then,
$$
\det (E_{p+1})=\det(E_{p+1}/E_p)\det (E_{p})=(w-e_i^tE_p^{-1}e_j)\det (E_{p})=-\det (E_{p})\neq 0,
$$
$$
\det(I_{p+1}+E_{p+1})=(1+w-e_i^t(I_p+E_p)^{-1}e_j)\det(I_p+E_{p})$$$$=
(1+w-e_i^t(E_p^{-1}-R)e_j)\det(I_p+E_{p})
=R(i,j)\det(I_p+E_{p})\neq 0,
$$
hence, $E_{p+1}, I_{p+1}+E_{p+1}, \in \Gl_{p+1}(\FF)$.

\hfill $\Box$
\begin{remark}
If $\mid\FF\mid \neq 2$, Lemma \ref{lemmaei+e} is straightforward,  and the result holds for $n\geq 1$. We can take, for example,   $E_n=cI_n$, with $c\in \FF$, $c\neq 0, -1$.
\end{remark}

\begin{proposition}\label{proprn}
Let $n\geq 2$ and
$A(s)=sI_n+A, B(s)=sI_n+B \in \FF[s]^{n\times n}$. Then there exists a matrix pencil
$P(s)\in \FF[s]^{n \times n}$ with $\rank (P(s))=n$ such that $A(s)+P(s)\se  B(s)$.
\end{proposition}

{\bf Proof.}
By Lemma \ref{lemmaei+e}, there exists $E_n\in \Gl_n(\FF)$ such that $I_n+E_n\in \Gl_n(\FF)$.

Let  $P_0=(I_n+E_n)B-A \in \FF^{n\times n}$
and  $P(s)=E_ns+P_0\in \FF[s]^{n\times n}$. Then  $\rank (P(s))=n$ and
$$
A(s)+P(s)=sI_n+A+E_ns+(I_n+E_n)B-A=(I_n+E_n)(sI_n+B)\se sI_n+B.
$$
\hfill $\Box$

When $r<n$, next lemma    allows us to take advantage of a solution to the bounded case and out of it  to built a solution for the fixed rank case. This is done in Proposition \ref{mainprop}.

\begin{lemma}\label{mainlemma}
Let $r_1, r, n$ be integers, $0\leq r_1<r< n$.
Let $I, J\subset \{1, \dots, n\}$ be such that $\mid I\mid=\mid J\mid  =r_1\geq 0$.  Then there exists a matrix $E\in \FF^{n \times n}$ satisfying that
$\rank (E)=r-r_1$, $I_n+E\in \Gl_n(\FF)$,  $E(I, :)=0$, and $E(:, J)=0$.
\end{lemma}

{\bf Proof.}
First, let us see that there exist sets
$$
R_1=\{i_1, \dots, i_{x'}\}, \quad
R_2=\{i_{x'+1}, \dots, i_{x'+a'}\},\quad S_2=\{i_{x'+a'+1}, \dots, i_{x'+2a'}\}
$$
($x'\geq 0$, $a'\geq 0$) with $i_k\neq i_\ell$ for  $k\neq \ell$, such that $R_1\dot {\cup} R_2 \subset I^c$, $R_1\dot{\cup} S_2 \subset J^c$,  $x'+a'=r-r_1$, and $x'\neq 1$.

Let $$X=I^c\cap J^c,\quad Y=I^c\setminus X,\quad Z=J^c\setminus X$$
and let $x=\mid  X\mid $, $a=n-r_1-x=\mid  Y\mid =\mid  Z\mid $.

\bigskip
\begin{itemize}
  \item
If  $a\geq r-r_1$, we put
$R_1=\emptyset$ and  choose $R_2\subseteq Y$, $S_2\subseteq Z$ such that $\mid  R_2\mid = \mid S_2 \mid =r-r_1$.
 In this case, $x'=0$, $a'=r-r_1$.
\item
If  $a<r-r_1$,  then $x=n-r_1-a>n-r> 0$. Therefore
$x\geq 2$.

\begin{itemize}
  \item
If  $(r-r_1)-a\geq 2$ we put $R_2=Y$, $S_2=Z$ and choose $R_1\subset X$ with  $\mid R_1\mid =r-r_1-a(<n-r_1-a=x)$.  In this case, $x'=r-r_1-a\geq 2$, $a'=a$.
  \item
    If  $(r-r_1)-a= 1$ and  $a\geq 1$,   we choose $R_1\subseteq X$ with  $\mid R_1\mid =2$ and $R_2\subset Y$, $S_2\subset Z$ with $\mid R_2\mid =\mid S_2\mid =r-r_1-2=a-1$.
    In this case, $x'= 2$, $a'=a-1$.
  \item
    If  $(r-r_1)-a= 1$ and $a=0$, then $r-r_1=1< x$. We can choose $i, j\in X$ such that  $i\neq j$. We put
    $R_1=\emptyset$, $R_2=\{i\}$, $S_2=\{j\}$.
    In this case, $x'= 0$, $a'=1$.    
\end{itemize}

\end{itemize}

We have that
$
R_1\dot{\cup}R_2\dot{\cup}S_2\subseteq \{1, \dots, n\}
$, hence $x'+2a'\leq n$.
Let us denote 
 $(R_2\cup R_1\cup S_2)^c=\{i_{x'+2a'+1}, \dots, i_n\}$.

We have obtained that  $x'= 0$ or  $x'\geq 2$. If $x'\geq 2$, by Lemma \ref{lemmaei+e}
there exists
$E_{x'}\in \Gl_{x'}(\FF)$
such that  $I_{x'}+E_{x'}\in \Gl_{x'}(\FF)$.

Let $\bar E\in \FF^{n\times n}$ be the matrix having
$$
\bar E(\{1, \dots,  x'\},\{1, \dots,  x'\} )=E_{x'}, \;
\bar E( \{x'+1, \dots,   x'+a'\}, \{x'+a'+1, \dots, x'+2a'\})=I_{a'},$$
and the rest of its entries equal to zero, i.e.
$$
\bar E=
\begin{bmatrix}
  E_{x'}&0&0&0\\
  0&0&I_{a'}&0\\
  0&0&0&0\\
  0&0&0&0\\
\end{bmatrix}\in \FF^{(x'+a'+a'+(n-x'-2a'))\times (x'+a'+a'+(n-x'-2a'))}.
$$
(If $x'=0$ or $a'=0$, the corresponding block vanishes).
Obviously,
$$\rank  (\bar E)=\rank  (E_{x'})+\rank (I_{a'})=x'+a'=r-r_1,$$
and
$$ I_n+\bar E=
\begin{bmatrix}
  I_{x'}+E_{x'}&0&0&0\\
  0&I_{a'}&I_{a'}&0\\
  0&0&I_{a'}&0\\
  0&0&0&I_{n-x'-2a'}
  \end{bmatrix}
\in \Gl_n(\FF).$$
Let $P$ be  the permutation matrix
$P=\begin{bmatrix}
e_{i_1}&\dots  &e_{i_n}
\end{bmatrix}
$, where $e_k$ denotes 
the $k$-th  column of  $I_n$. Then,  $Pe_k=e_{i_k}$ for $1\leq k \leq n$; equivalently,  $P^te_{i_k}=e_k$, and $e_{i_k}^tP=e_k^t$.

Let
$E=P\bar E P^t$.
Then,
$\rank (E)=\rank (\bar E)=r-r_1$, $I_n+ E=PP^t+P\bar E P^t=P(I_n+\bar E)P^t\in \Gl_n(\FF),$
$$
E((R_1\cup R_2)^c, :))=E(\{i_{x'+a'+1}, \dots, i_n\}, :)=\begin{bmatrix}
e_{i_{x'+a'+1}}^t\\\vdots \\e_{i_n}^t
\end{bmatrix}P\bar E P^t=\begin{bmatrix}
e_{x'+a'+1}^t\\\vdots \\e_n^t
\end{bmatrix}\bar E P^t$$$$=\bar E(\{x'+a'+1, \dots, n\}, :)P^t=0,
$$
and
$$
E(:, (R_1\cup S_2)^c)=
E(:,\{i_{x'+1}, \dots, i_{x'+a'}\}\cup\{i_{x'+2a'+1}, \dots, i_{n}\})$$$$=
P\bar E P^t
\begin{bmatrix}
e_{i_{x'+1}}&\dots &e_{i_{x'+a'}}&
e_{i_{x'+2a'+1}}&\dots &e_{i_n}
\end{bmatrix}$$$$=P\bar E \begin{bmatrix}
e_{x'+1}&\dots &e_{x'+a'}&
e_{x'+2a'+1}&\dots &e_n
\end{bmatrix}$$$$=P\bar E(:, \{x'+1, \dots, x'+a'\}\cup \{x'+2a'+1, \dots, n\})=0.
$$
Since  $I\subseteq (R_1\cup R_2)^c$ and $J\subseteq (R_1\cup S_2)^c$, it results  that
$
E(I, :)=0$ and  $E(:, J)=0.
$

\hfill $\Box$

\begin{proposition}\label{mainprop}
Let $n\geq 2$ and
$A(s)=sI_n+A \in \FF[s]^{n\times n}$. Let $P\in \FF^{n \times n}$ be a matrix  such that $\rank (P)=r_1$
and
let $r$ be an integer, $r_1<r<n$. Then there exists a matrix pencil
$P(s)\in \FF[s]^{n \times n}$ with $\rank (P(s))=r$ such that $A(s)+P(s)\se A(s)+P$.
\end{proposition}

{\bf Proof.}
Since  $\rank (P)=r_1$, there exist  $I, J\subset \{1, \dots, n\}$ such that
$\mid I\mid =  \mid J\mid =r_1$ and  $\det (P(I, J))\neq 0$ (if $r_1=0$, then
$I=J=\emptyset$).
By Lemma \ref{mainlemma}, there exists a matrix $E\in \FF^{n \times n}$ such that
$\rank (E)=r-r_1$, $I_n+E\in \Gl_n(\FF)$,  $E(I, :)=0$, and $E(:, J)=0$.

Let  $Q=I_n+E$. Then
$$
sI_n+A+P\se Q(sI_n+A+P)=sI_n+A+P+E(sI_n+A+P)=A(s)+P(s),
$$
where $P(s)=P+E(sI_n+A+P)$.

Let us see that  $\rank (P(s))=r$.
On one hand,
$$\rank (P(s))\leq \rank (P)+ \rank (E(sI_n+A+P))\leq \rank (P)+ \rank (E)=r_1+r-r_1=r.$$
On the other one,
$$
P(s)(I, :)=P(I, :)+E(I,:)(sI_n+A+P)=P(I,:).
$$
Therefore,
$$\det (P(s)(I, J))=\det (P(I, J))\neq 0,$$
and
$$P(s)/P(s)(I, J)=P(s)(I^c, J^c)-P(s)(I^c, J)  P(I, J)^{-1}P(I, J^c).$$
As
$$
P(s)(I^c, J^c)=P(I^c,J^c)+sE(I^c,J^c)+(E(A+P))(I^c, J^c),
$$
and
$$
P(s)(I^c,J)=P(I^c,J)+sE(I^c,J)+(E(A+P))(I^c, J)$$$$=P(I^c,J)+(E(A+P))(I^c, J)\in \FF^{(n-r_1)\times r_1},
$$
we can write
$P(s)/P(s)(I, J)=sE(I^c, J^c)+P_0,$ with $P_0\in \FF^{(n-r_1)\times (n-r_1)}$, from where
$$\rank (P(s)/P(s)(I, J))\geq \rank (E(I^c, J^c))=\rank (E)=r-r_1.$$ Hence,
$$
\rank (P(s))=\rank (P(I, J))+\rank (P(s)/P(s)(I, J))\geq r.
$$

\hfill $\Box$

\begin{theorem}\label{propsI+A}
Let $n\geq 2$ and  $A(s)=sI_n+A, B(s)=sI_n+B\in \FF[s]^{n \times n}$.
Let $\phi_1(s,t)\mid \dots\mid \phi_n(s, t)$ and $\psi_1(s,t)\mid \dots\mid \psi_n(s, t)$ be the homogeneous invariant factors of $A(s)$ and $B(s)$, respectively. Let $r$ be a nonnegative integer,
$r\leq n$. If (\ref{interlacinghomogr1}) is satisfied, then there exists
a matrix pencil
$P(s)\in \FF[s]^{n \times n}$ such that $\rank (P(s))= r$ and
$A(s)+P(s)\se B(s)$.
\end{theorem}

{\bf Proof.}
If $r=n$, we apply Proposition \ref{proprn}.

If  $r<n$,  let $\alpha_i(s)=\phi_i(s, 1)$ and $\beta_i(s)=\psi_i(s, 1)$,
 $1\leq i \leq n$,
be
the invariant factors of $A(s)$  and $B(s)$, respectively.
Then, conditions
 (\ref{interlacinghomogr1}) imply conditions (\ref{eqintif}).
By Corollary \ref{corless}, there exists $P\in \FF^{n \times n}$ such that $\rank (P)\leq r$ and
$A(s)+P\se B(s)$. If $\rank (P)< r$, we apply Proposition \ref{mainprop}.

\hfill $\Box$

We show next an example of pencils $A(s)$ and $B(s)$ such that $B(s)$ cannot be obtained by a constant perturbation of rank $2$  of $A(s)$, but it does result as a pencil perturbation of rank $2$ of the pencil $A(s)$.

\begin{example}\label{exp1p}
Let $\FF$ be an arbitrary field and $r=2$,
$$A(s)=\begin{bmatrix}s-1& 0&0\\0&s-1&0\\0&0&s-1 \end{bmatrix}, \quad
  B(s)=\begin{bmatrix}s-1& 0&0\\0&s-1&0\\0&0&s \end{bmatrix}.$$
The homogeneous invariant factors of $A(s)$ and $B(s)$ are
  $\phi_1(s, t)=\phi_2(s, t)=\phi_3(s, t)=(s-t)$ and
  $\psi_1(s, t)=1, \psi_2(s, t)=(s-t),  \psi_3(s, t)=s(s-t)$, respectively.
We have that
  $$
 \phi_{i-2}(s, t)\mid \psi_i(s, t)\mid\phi_{i+2}(s, t), \quad 1\leq i\leq 3.
  $$
 Therefore, 
 $$
 \phi_{i-2}(s, 1)\mid \psi_i(s, 1)\mid\phi_{i+2}(s, 1), \quad 1\leq i\leq 3,
  $$
hence,   by Corollary \ref{corless}, there exists a matrix $P\in \FF^{3\times 3}$ such that $\rank P\leq 2$ and
 $A(s)+P \se B(s)$.
 In fact, taking $P= \begin{bmatrix}0& 0&0\\0&0&0\\0&0&1 \end{bmatrix}$, we have that $\rank P =1$ and $A(s)+P = B(s)$.

 Observe that
 $$
\min\{\rank A(\lambda)+ \rank B(\lambda): \ \lambda \in \bar{\FF}\}=1.
$$
By Proposition \ref{propSi2}, this means that there is no  $P\in \bar{\FF}^{3\times 3}$ such that $\rank P= 2$ and $A(s)+P \se B(s)$.

Let $Q= \begin{bmatrix}1& 1&0\\0&1&0\\0&0&1 \end{bmatrix}\in \Gl_3(\FF)$.
Then
$$
B(s) \se QB(s)=Q(A(s)+P)=\begin{bmatrix}s-1& s-1&0\\0&s-1&0\\0&0&s \end{bmatrix}=A(s)+P(s),
$$
where $P(s)=\begin{bmatrix}0& s-1&0\\0&0&0\\0&0&1 \end{bmatrix}\in \FF[s]^{3 \times 3}$, and $\rank P(s)=2$.
  \end{example}

\begin{corollary}\label{cornotid}
Let $n\geq 2$ and let $A(s)=A_0+sA_1, B(s)=B_0+sB_1\in \FF[s]^{n \times n}$ be such that $\det (A_1)\neq 0$ and $\det (B_1) \neq 0$.
Let $\phi_1(s,t)\mid \dots\mid \phi_n(s, t)$ and $\psi_1(s,t)\mid \dots\mid \psi_n(s, t)$ be the homogeneous invariant factors of $A(s)$ and $B(s)$, respectively. Let $r$ be a nonnegative integer,
$r\leq n$. If (\ref{interlacinghomogr1}) is satisfied, then there exists
a matrix pencil
$P(s)\in \FF[s]^{n \times n}$ such that $\rank (P(s))= r$ and
$A(s)+P(s)\se B(s)$.
\end{corollary}

{\bf Proof.}
We have that $A(s)\se  A_1^{-1}A_0+sI_n$ and $B(s)\se B_1^{-1}B_0+sI_n$.
Hence, the homogeneous invariant factors of $sI_n+A_1^{-1}A_0$ and $sI_n+B_1^{-1}B_0 $ are $\phi_1(s,t)\mid \dots\mid \phi_n(s, t)$ and $\psi_1(s,t)\mid \dots\mid \psi_n(s, t)$, respectively.
By Theorem \ref{propsI+A}, there exists a matrix pencil $P'(s)\in \FF[s]^{n \times n}$ such that
$\rank (P'(s))=r$ and  $sI_n+A_1^{-1}A_0+P'(s)\se sI_n+B_1^{-1}B_0\se B(s)$.
Let $P(s)=A_1P'(s)$. Then, $\rank (P(s))=\rank (P'(s))=r$ and
$$
A(s)+P(s)=A_1(sI_n+A_1^{-1}A_0+P'(s))\se sI_n+A_1^{-1}A_0+P'(s)
\se B(s).
$$
\hfill $\Box$

\subsection{General case}
\label{subsecgeneral}

We analyze first the case $n=1$.

\begin{theorem} \label{theon1}
Let $a(s)=a_0+sa_1, b(s)=b_0+sb_1\in \FF[s]$ be such that $a(s)\neq 0$ and  $b(s)\neq 0$.
Let $\phi_1(s,t)$ and $\psi_1(s,t)$ be the homogeneous invariant factors of
$a(s)$ and $b(s)$, respectively. Let $r$ be an integer, $0\leq r\leq 1$.
\begin{enumerate}
\item If  $\mid\FF\mid >2$ or  $r=0$, then there exists $p(s)=p_0+sp_1\in \FF[s]$ such that $\rank (p(s))=r$ and $a(s)+p(s)\se b(s)$ if and only if (\ref{interlacinghomogr1}) holds.
\item
If $\mid\FF\mid =2$ and  $r=1$, then there exists $p(s)\in \FF[s]$ such that $\rank (p(s))=1$ and $a(s)+p(s)\se b(s)$ if and only if $a(s)\neq b(s)$.
\end{enumerate}
\end{theorem}

{\bf Proof.}
\begin{enumerate}
\item
The necessity is an immediate consequence  of Proposition \ref{propless}. Let us prove the sufficiency.

Since $n=1$,  conditions (\ref{interlacinghomogr1}) reduce to
\begin{equation}\label{eqn1}
  \psi_{1-r}(s,t)\mid\phi_1(s,t)\mid\psi_{1+r}(s,t).
  \end{equation}
\begin{itemize}
\item
  If $r=0$, then (\ref{eqn1}) implies  $\psi_{1}(s)=\phi_1(s)$, hence $b(s)\se a(s)=a(s)+0$.
\item
  If $r=1$, then (\ref{eqn1}) is trivially satisfied for any $a(s), b(s)$. 
  As  $\mid\FF\mid >2$, there exists $c\in \FF\setminus\{ 0\}$ such that $a(s)\neq cb(s)$. Taking  $p(s)=cb(s)-a(s)$, the sufficiency is proven. 
    
\end{itemize}

\item  It is enough to observe that if $\mid\FF\mid =2$, there exists $p(s)\in \FF[s]$ such that  $a(s)+p(s)\se b(s)$ if and only if $a(s)+p(s)= b(s)$.
\end{enumerate}

\hfill $\Box$

Next theorem is our main result.

\begin{theorem}
  \label{maintheorem}
Let $n\geq 2$. Let $A(s)=sA_1+A_0, B(s)=sB_1+B_0\in \FF[s]^{n \times n}$ be regular matrix pencils.
Let $\phi_1(s,t)\mid \dots\mid \phi_n(s, t)$ and $\psi_1(s,t)\mid \dots\mid \psi_n(s, t)$ be the homogeneous invariant factors of $A(s)$ and $B(s)$, respectively, and assume that  $ \FF\cup \{\infty\}\not \subseteq \Lambda(A(s))\cup \Lambda(B(s))  $. Let $r$ be a nonnegative integer, $r\leq n$.
There exists a matrix pencil  $P(s)\in \FF[s]^{n \times n}$ such that $\rank (P(s))=r$ and $A(s)+P(s)\se B(s)$ if and only if
(\ref{interlacinghomogr1}) holds.

\end{theorem}

{\bf Proof.}
 The necessity is an immediate consequence  of Proposition \ref{propless}.

Assume that (\ref{interlacinghomogr1}) holds.
As $ \FF\cup \{\infty\}\not \subseteq \Lambda(A(s))\cup \Lambda(B(s)) $, there exists $c\in  \FF\cup \{\infty\}$ such that $c\not \in \Lambda(A(s))\cup \Lambda(B(s))$.

If $c=\infty$, we apply Corollary \ref{cornotid}.

If $c\neq \infty$, take
$$
X=\begin{bmatrix}
c&1\\
1&0
\end{bmatrix},
$$
and
$$A'(s)=P_X(sA_1+A_0)=s(cA_1+A_0)+A_1=sA_1'+A_0',$$
$$B'(s)=P_X(sB_1+B_0)=s(cB_1+B_0)+B_1=sB_1'+B_0'.$$
Then, $\det (A'_1)\neq 0$, $\det (B'_1)\neq 0$.

 Let $\phi'_1(s, t), \dots, \phi'_n(s, t)$ and
$\psi'_1(s, t), \dots, \psi'_n(s, t)$ be the homogeneous invariant factors of $A'(s)$ and $B'(s)$, respectively.
By Lemma \ref{lemmaCasi10} and Remark \ref{remrank},
$$
\phi'_i(s, t)=c_i\Pi_X(\phi_i)(s,t), \quad \psi'_i(s, t)=d_i\Pi_X(\psi_i)(s,t),\quad 1\leq i \leq n,
$$
where $0\neq c_i\in \FF$, $0\neq d_i\in \FF$, $1\leq i \leq n$.
Applying Lemma \ref{lemmaCasi7}, from (\ref{interlacinghomogr1}) we obtain
$$
\phi'_{i-r}(s, t)\mid \psi'_i(s, t)\mid\phi'_{i+r}(s, t), \quad 1\leq i \leq n.
$$

By Corollary \ref{cornotid}, there exists a matrix pencil $P'(s)=sP_1'+P_0'\in \FF[s]^{n \times n}$ such that $\rank (P'(s))=r$ and
$$A'(s)+P'(s)\se B'(s).$$
Then, by Lemmas \ref{lemmaCasi6} and \ref{lemmaCasi9}, 
$$\left(P_X\right)^{-1}(A'(s))+\left(P_X\right)^{-1}(P'(s))=\left(P_X\right)^{-1}(A'(s)+P'(s)) \se \left(P_X\right)^{-1}(sB_1'+B_0').$$
Taking  $P(s)=\left(P_X\right)^{-1}(P'(s))=P_{X^{-1}}(P'(s))$, we obtain that 
$$
A(s)+P(s)\se  B(s),
$$
and by Remark \ref{remrank}, $\rank (P(s))=r$.

\hfill $\Box$

\begin{corollary}\label{corrminfixed2}
Let $n\geq 2$. Let $A(s)=sA_1+A_0, B(s)=sB_1+B_0\in \FF[s]^{n \times n}$ be regular matrix pencils.
Let $\phi_1(s,t)\mid \dots\mid \phi_n(s, t)$ and $\psi_1(s,t)\mid \dots\mid \psi_n(s, t)$ be the homogeneous invariant factors of $A(s)$ and $B(s)$, respectively, and assume that  $ \FF\cup \{\infty\}\not \subseteq \Lambda(A(s))\cup \Lambda(B(s))$. 
Let
$$
r_0=\min\{r\geq 0 \; : \;  \phi_{i-r}(s,t)\mid\psi_{i}(s,t)\mid\phi_{i+r}(s,t), \quad 1\leq i \leq n\}.
$$
Then
there exists a matrix pencil  $P(s)\in \FF[s]^{n \times n}$ with $\rank (P(s))=r$ and such that $A(s)+P(s)\se B(s)$ if and only if 
$r_0\leq r\leq n$.
\end{corollary}

{\bf Proof.}
It is straigthtforward that $r\geq r_0$ if and only if conditions (\ref{interlacinghomogr1}) hold. 

\hfill $\Box$

\begin{example}
Let $\FF$ be an arbitrary field. Let $A(s), B(s)\in \FF[s]^{5\times 5}$ be regular matrix pencils with homogeneous invariant factors 
$$\phi_1(s,t)=\phi_2(s,t)=1, \quad \phi_3(s,t)=t, \quad \phi_4(s,t)=\phi_5(s,t)=t^2,$$
$$\psi_1(s,t)=\psi_2(s,t)=1, \quad \psi_3(s,t)=\psi_4(s,t)=s-t,\quad
\psi_5(s,t)=(s-t)^3,$$
respectively. Then
$$
\Lambda(A(s))=\{\infty\}, \quad \Lambda(B(s))=\{1\}, \quad 0\not \in \Lambda(A(s))\cup\Lambda(B(s)),
$$
and
$$
r_0=\min\{r\geq 0 \; : \;  \phi_{i-r}(s,t)\mid\psi_{i}(s,t)\mid\phi_{i+r}(s,t), \quad 1\leq i \leq 5\}=3.
$$
Hence, for $3\leq r \leq 5$
there exist matrix pencils $P_r(s)\in \FF[s]^{5\times 5}$ with  
$\rank(P_r(s))=r$  such that $A(s)+P_r(s)\se B(s)$.

Moreover, there is not any pencil $P(s)$  with  $\rank(P(s))\leq 2$  such that $A(s)+P(s)\se B(s)$.
\end{example}

The characterization of the solution given in Theorem \ref{maintheorem} can be stated in terms of the partial multiplicities of the elements of $ \Lambda(A(s))\cup \Lambda(B(s))$  
(see \cite[Corollary 4.5]{BaRo18} for an analogous result when $\rank (P(s))\leq r$; see also \cite[Proposition 4.2]{GeTr17}  for $r=1$).

\begin{corollary}\label{cormainmult}
Let $n\geq 2$. Let $A(s), B(s)\in \FF[s]^{n \times n}$ be regular matrix pencils.
Assume that  $ \FF\cup \{\infty\}\not \subseteq \Lambda(A(s))\cup \Lambda(B(s))  $.
Let $r$ be a nonnegative integer, $r\leq n$.
There exists a matrix pencil  $P(s)\in \FF[s]^{n \times n}$ such that $\rank (P(s))= r$ and $A(s)+P(s)\se B(s)$ if and only if
\begin{equation}\label{eqnecmult}
m_{i-r}(\lambda, A(s))\leq m_i(\lambda, B(s))\leq m_{i+r}(\lambda, A(s)),
\, 1\leq i \leq n,  \quad \lambda \in \overline{\FF}\cup\{\infty\}.
\end{equation}

\end{corollary}

As pointed out in \cite[Remark 4.15]{BaRo18}, if $\#\FF> 2n $, 
 the condition $\FF\cup \{\infty\} \not \subseteq \Lambda(A(s))\cup \Lambda(B(s))$ is automatically satisfied. In  the case that    $\#\FF\leq 2n $, Theorem \ref{maintheorem} can still be applied if there exists an element $c\in \FF\cup \{\infty\}$ which is neither an eigenvalue of  $A(s)$ nor of $B(s)$.

Moreover, we show in Corollary \ref{coreqmult} that the condition $ \FF\cup \{\infty\} \not \subseteq \Lambda(A(s))\cup \Lambda(B(s))$ is not always necessary.

 \begin{corollary}\label{coreqmult} 
Let $A(s), B(s)\in \FF[s]^{n \times n}$ be regular matrix pencils.
Let $\phi_1(s,t)\mid \dots\mid \phi_n(s, t)$ and $\psi_1(s,t)\mid \dots\mid \psi_n(s, t)$ be the homogeneous invariant factors of 
$A(s)$ and $B(s)$, respectively, and assume that  
for some $\lambda_0\in \FF\cup \{\infty\}$,
$$
m_i(\lambda_0, A(s))=m_i(\lambda_0, B(s)), \quad 1\leq i \leq n.
$$
Let $r$ be a nonnegative integer, $r\leq n$. 
There exists a matrix pencil  $P(s)\in \FF[s]^{n \times n}$ such that $\rank (P(s))=r$ and $A(s)+P(s)\se B(s)$ if and only if
(\ref{interlacinghomogr1}) holds.
\end{corollary}

{\bf Proof.}
Analogous to the proof  of \cite[Theorem 4.17]{BaRo18}. \hfill$\Box$

\begin{example}\label{exeqmult}
Let $\FF=\ZZ_2$, $r=2$,
$$\hat A(s)=\begin{bmatrix}1&0&0&0\\0&s-1& 0&0\\0&0&s-1&0\\0&0&0&s-1 \end{bmatrix}, \quad
  \hat B(s)=\begin{bmatrix}1&0&0&0\\0&s-1& 0&0\\0&0&s-1&0\\0&0&0&s \end{bmatrix}.$$
The homogeneous invariant factors of $\hat A(s)$ and $\hat B(s)$ are
  $\phi_1(s, t)=1, \phi_2(s, t)=\phi_3(s, t)=(s-t),\phi_4(s, t)=t(s-t) $ and
  $\psi_1(s, t)=\psi_2(s, t)=1, \psi_3(s, t)=(s-t),  \psi_4(s, t)=ts(s-t)$, respectively.
 Then
  $$
 \phi_{i-2}(s, t)\mid \psi_i(s, t)\mid\phi_{i+2}(s, t), \quad 1\leq i\leq 4,
  $$
$$
\Lambda(\hat A(s))=\{1, \infty\}, \quad \Lambda(\hat B(s))=\{0, 1, \infty\},
$$
and  $ \FF\cup \{\infty\}= \Lambda(\hat A(s))\cup \Lambda(\hat B(s))=\{0,1, \infty\}$.
But
$$
(m_1(\infty,\hat  A(s)), \dots, m_4(\infty,\hat  A(s)))=(m_1(\infty, \hat B(s)), \dots, m_4(\infty, \hat B(s)))
=(0,0,0,1).
$$
We have that
$$
\hat A(s)=\begin{bmatrix}1&0\\0&A(s)\end{bmatrix}, \quad
\hat B(s)=\begin{bmatrix}1&0\\0&B(s)\end{bmatrix},
$$
where $A(s)$ and $B(s)$ are the pencils of Example \ref{exp1p} and we have seen that there exists a matrix pencil $P(s)\in \FF[s]^{3\times 3}$  such that
$\rank P(s)= 2$ and 
$A(s)+P(s)\se B(s)$. Taking 
$\hat P(s)= \begin{bmatrix}0&0\\0&P(s)\end{bmatrix}\in \FF[s]^{(1+3)\times (1+3)}$, we have  that
$\hat A(s)+\hat P(s)\se \hat B(s)$ and $\rank \hat P(s)=2$.

  \end{example}

\section{Eigenvalue placement for regular matrix pencils under fixed rank perturbations}
\label{secplacement}
In this section we give a solution to Problem \ref{prfixdet}.

\medskip
Recall that if $\Gamma(s, t)$ is an homogeneous polynomial, 
$$\Lambda(\Gamma(s,t)):=\{\lambda \in \overline{\FF}\cup \{\infty\}\; : \; 
\Gamma(\lambda,1)=0\},$$
where $\Gamma(\infty, 1):=\Gamma(1, 0)$.

The following theorem is the main result  in this section.
The proof is similar to that of Theorem 5.1 of \cite{BaRo18}.

\begin{theorem}\label{thomainplacement}
Let $n\geq 2$.
Let $A(s)\in \FF[s]^{n \times n}$ be a regular matrix pencil and $\phi_1(s,t)\mid \dots\mid \phi_n(s,t)$ be its homogeneous  invariant factors. Let $\Psi(s,t)\in \FF[s,t]$ be a nonzero homogeneous polynomial, monic with respect to $s$, and such that $\deg (\Psi(s, t))= n$.
Assume that  $ \FF\cup \{\infty\}\not \subseteq \Lambda(A(s))\cup 
\Lambda(\Psi(s,t))  $.
Let $r$ be  a nonnegative integer, $r\leq n$.
There exists a matrix pencil  $P(s)\in \FF[s]^{n \times n}$ with $\rank (P(s))=r$ 
such that if $C(s, t)$ is the homogeneous pencil associated to $A(s)+P(s)$, then $\det (C(s,t))=k\Psi(s,t)$ with $0\neq k \in \FF$ if and only if
\begin{equation}\label{desphi}
\phi_1(s,t)\dots\phi_{n-r}(s,t)\mid \Psi(s,t).
\end{equation}

\end{theorem}

{\bf Proof.}
{\it Necessity.}
Let $C(s)=A(s)+P(s)$ and let 
$\psi_1(s,t)\mid \dots\mid \psi_n(s,t)$ be its homogeneous invariant factors.
Taking  $\Psi(s,t)=\psi_1(s,t) \dots \psi_n(s,t)$,  from Theorem \ref{maintheorem} condition  (\ref{desphi}) is satisfied.

\medskip

{\it Sufficiency.}
Assume that (\ref{desphi}) holds.
Then, there exists an homogeneous polynomial 
$\gamma(s, t)\in \FF[s, t]$ such that
$$\Psi(s,t)=\phi_1(s,t)\dots\phi_{n-r}(s,t)\gamma(s, t).$$
We define
    $$
\psi_i(s,t):=\phi_{i-r}(s,t), \; 1\leq i \leq n-1, \quad \psi_n(s,t):=
\phi_{n-r}(s, t)\gamma(s, t),
$$
then
$$\psi_1(s,t) \mid \dots \mid \psi_n(s,t)   \ \text{ and } \
\sum_{i=1}^n\deg(\psi_i(s, t))=\deg (\Psi(s, t))=n.
$$

Let  $B(s)$ be a pencil with homogeneous invariant factors
 $\psi_1(s,t) \mid \dots \mid \psi_n(s,t)$.
Then,  $B(s)$ is regular and condition  (\ref{interlacinghomogr1}) is satisfied.
By Theorem \ref{maintheorem},
there exists a pencil  $P(s)\in \FF[s]^{n \times n}$ such that  $\rank (P(s))= r$ and $A(s)+P(s)\se B(s)$.
 Let $B(s,t)$  be the homogeneous pencil associated to $B(s)$.
Then there exist $0\neq k_1, k_2 \in \FF$ such that
$$
\det (C(s,t))=k_1\det (B(s,t)) =k_1k_2\psi_1(s,t)\dots \psi_n(s,t)$$$$=k_1k_2\phi_1(s,t)\dots\phi_{n-r}(s,t)\gamma(s,t)=k_1k_2\Psi(s,t),\quad 0\neq k_1k_2\in \FF.
$$
\hfill $\Box$

Notice that Theorem \ref{thomainplacement} gives us a solution to Problem \ref{prfixdet} as we see in the following corollary (compare it with Theorem 5.4 in \cite{BaRo18}).

\begin{corollary}\label{corplacement}
Let $n\geq 2$.
Let $A(s)\in \FF[s]^{n \times n}$ be a regular matrix pencil and $\alpha_1(s)\mid \dots\mid \alpha_n(s)$ be its   invariant factors. Let $q(s)\in \FF[s]$ be a nonzero monic polynomial with $\deg (q(s))\leq n$.
Assume that  $ \FF\cup \{\infty\}\not \subseteq \Lambda(A(s))\cup \Lambda^n(q(s))  $.
Let $r$ be  a nonnegative integer, $r\leq n$.
There exists a matrix pencil  $P(s)\in \FF[s]^{n \times n}$ such that $\rank (P(s))= r$ and $\det(A(s)+P(s))=kq(s)$ with $0\neq k \in \FF$ if and only if
\begin{equation}\label{eqdetpol}
\begin{array}{l}
\alpha_1(s)\dots\alpha_{n-r}(s)\mid q(s),\\
\sum_{i=1}^{n-r}m_i(\infty, A(s)) \leq n-\deg(q(s)).
\end{array}
\end{equation}

\end{corollary}

{\bf Proof.}
Let  $\phi_1(s,t)\mid \dots\mid \phi_n(s,t)$ be the homogeneous   invariant factors of $A(s)$ and
let $\Psi(s,t)=t^{n}q(\frac{s}{t})$.
Then $\Psi(s,t)\in \FF[s,t]$ is a nonzero homogeneous polynomial,  $\deg (\Psi(s, t))= n$ and
$ \FF\cup \{\infty\}\not \subseteq \Lambda(A(s))\cup 
\Lambda(\Psi(s,t))  $.
Take $\delta(s)=\alpha_1(s)\dots\alpha_{n-r}(s)$. Then
$$
\phi_1(s,t)\dots\phi_{n-r}(s,t)=t^{\sum_{i=1}^{n-r}m_i(\infty, A(s))}t^{\deg(\delta)}\delta(\frac{s}{t}).
$$
Hence, (\ref{desphi}) is equivalent to (\ref{eqdetpol}).

Assume that there exists a matrix pencil  $P(s)\in \FF[s]^{n \times n}$ such that $\rank (P(s))= r$ and $\det(A(s)+P(s))=kq(s)$ with   $0\neq k \in \FF$. Let
$C(s, t)$ be the homogeneous pencil associated to $C(s)=A(s)+P(s)$.
Then, $\deg(\det (C(s, t)))=n$ and $\det (C(s, 1))=\det(A(s)+P(s))=kq(s)$, from where 
$\det (C(s, t))=t^{n-\deg(q)}t^{\deg(q)}kq(\frac{s}{t})=k\Psi(s,t)$. By Theorem \ref{thomainplacement}, (\ref{desphi}) (equivalently, 
(\ref{eqdetpol})) holds.

Conversely, assume that (\ref{eqdetpol}) (equivalently, 
(\ref{desphi})) holds. Then by Theorem \ref{thomainplacement}, there exists a matrix pencil  $P(s)\in \FF[s]^{n \times n}$ with $\rank (P(s))=r$ and 
such that if $C(s, t)$ is the homogeneous pencil associated to $C(s)=A(s)+P(s)$, then $\det (C(s,t))=k\Psi(s,t)$ with $0\neq k \in \FF$.
Therefore, $\det (A(s)+P(s))=\det (C(s,1)) =kq(s)$.

\hfill $\Box$

\begin{remark}\label{rempln1}
If $n=1$, given pencils $a(s), q(s), p(s)\in \FF[s]$, then $\det(a(s)+p(s))=kq(s)$ with $0\neq k\in \FF$  if and only if  $a(s)+p(s)\se q(s)$, i.e.  Problem \ref{prfixdet} is the same as  Problem \ref{fixproblem}. 
\end{remark}

\section{Conclusions}
\label{conclusions}

Given a regular matrix pencil, we have completely characterized the
Weiestrass structure of a regular pencil obtained by a perturbation of  fixed rank. The characterization is stated in terms of interlacing conditions between the homogeneous invariant factors of the original an the perturbed pencils, except in a very particular case. This work completes  the research carried out in \cite{BaRo18}, where the same type of problem was solved  in the case that the perturbed pencil was of bounded rank. Surprisingly enough, both solutions are characterized in terms   of the same interlacing conditions.

The necessity of the conditions holds over  arbitrary fields and the sufficiency over fields with sufficient number of elements.

It is remarkable the fact that the characterization of the fix rank perturbation of a pencil of the form $A(s)=sI-A$ requires an extra condition when the perturbation is performed by a constant matrix (\cite{Silva88_1}), and that that extra condition disappears when the fixed rank perturbation is allowed to be a pencil of degree one.

We also solve an eigenvalue placement problem characterizing the assignment of the determinant  to a regular matrix pencil obtained by a fixed rank pencil perturbation   of another regular one.

\bibliographystyle{acm}
\bibliography{referencesreg}
\end{document}